\input amstex
 \input epsf
\magnification=\magstep1 
\baselineskip=13pt
\documentstyle{amsppt}
\def\SS{\Cal S}
\def\UU{\Cal U}
\vsize=8.7truein \CenteredTagsOnSplits \NoRunningHeads

\topmatter
\title More on zeros and approximation of the Ising partition function  \endtitle 
\author Alexander Barvinok  and Nicholas Barvinok \endauthor
\address Department of Mathematics, University of Michigan, Ann Arbor,
MI 48109-1043, USA \endaddress
\email barvinok$\@$umich.edu  \endemail
\thanks  The research of the first author was partially supported by NSF Grant DMS 1855428.
\endthanks 
\address School of Mathematics, Georgia Institute of Technology, Atlanta, GA 30332, USA \endaddress
\email nbarvinok3$\@$gatech.edu  \endemail
\date May 10, 2021 \enddate
\abstract We consider the problem of computing the partition function $\sum_x e^{f(x)}$, where $f: \{-1, 1\}^n \longrightarrow {\Bbb R}$ is a quadratic or cubic polynomial on the Boolean cube $\{-1, 1\}^n$. 
In the case of a quadratic polynomial $f$, we show that the partition function can be approximated within relative error $0 < \epsilon < 1$ in quasi-polynomial $n^{O(\ln n - \ln \epsilon)}$ time if the Lipschitz constant of the non-linear part of $f$ with respect to the $\ell^1$ metric on the Boolean cube does not exceed $1-\delta$, for any $\delta >0$, fixed in advance. For a cubic polynomial $f$, we get the same result under a somewhat stronger condition.
 We apply the method of polynomial interpolation, for which we prove that $\sum_x e^{\tilde{f}(x)} \ne 0$ for complex-valued polynomials $\tilde{f}$ in a neighborhood of a real-valued $f$ satisfying the above mentioned conditions. The bounds are asymptotically optimal. Results on the zero-free region are interpreted as the absence of a phase transition in the Lee - Yang sense in the corresponding Ising model.  The novel feature of the bounds is that they control the total interaction of each vertex but not every single interaction of sets of vertices.
\endabstract
\endtopmatter

\document

\head 1. Introduction and main results \endhead

The Ising model is one of the oldest, most famous and most studied models in statistical physics, see \cite{FV18} for a thorough introduction, description, results and references. In this paper, we look at the computational complexity and  complex zeros of the partition function in the Ising model. This is a classical and also currently very active area of research, see \cite{B+19}, \cite{G+19}, \cite{J+19}, \cite{L+12}, \cite{L+16}, \cite{L+19a}, \cite{L+19b}, \cite{PR20},
\cite{S+14} and \cite{Z+11}  for some recent results.

Formally, the partition function we work with is described as follows. Let $\{-1, 1\}^n$ be the $n$-dimensional Boolean cube of all $n$-vectors $x=\left(\xi_1, \ldots, \xi_n\right)$, where $\xi_i=\pm 1$ for 
$i=1, \ldots, n$ and let $f: \{-1, 1\}^n \longrightarrow {\Bbb C}$ be a function. We define the {\it partition function} of $f$ as 
$$\SS\left(e^f\right)=\sum_{x \in \{-1, 1\}^n} e^{f(x)}.$$ In this paper, we consider the cases of a quadratic 
$$f(x)=\sum_{1 \leq i < j \leq n} a_{ij} \xi_i \xi_j + \sum_{i=1}^n b_i \xi_i \quad \text{for} \quad x=\left(\xi_1, \ldots, \xi_n\right)$$
or cubic
$$f(x)=\sum_{1 \leq i < j < k \leq n} c_{ijk} \xi_i \xi_j \xi_k + \sum_{1 \leq i < j \leq n} a_{ij} \xi_i \xi_j + \sum_{i=1}^n b_i \xi_i \quad \text{for} \quad x=\left(\xi_1, \ldots, \xi_n\right)$$
polynomial $f$ and we are interested in computing (approximating) $\SS\left(e^f\right)$ efficiently.

The case of a quadratic polynomial $f$ is the best studied.
Often, there is an underlying graph $G$ with vertices numbered $1, \ldots, n$ and the property that $a_{ij} \ne 0$ if and only if $i$ and $j$ span an edge of $G$. In this case, the variable $\xi_i=\pm 1$ is interpreted as the {\it spin} of a vertex $i$ and $-f(x)$ as the {\it energy} of the {\it configuration} $x=\left(\xi_1, \ldots, \xi_n\right)$. The coefficients $a_{ij}$ describe the {\it interactions} of vertices $i$ and $j$: the interaction is {\it ferromagnetic} if $a_{ij} > 0$ and {\it antiferromagnetic} if $a_{ij} < 0$. The coefficients $b_i$ describe the {\it external field}, see \cite{FV18} for a thorough discussion. 

Via the method of polynomial interpolation \cite{Ba16}, we relate the computational complexity to the absence of complex zeros of $\SS\left(e^f\right)$ in the vicinity of real coefficients $a_{ij}$, $b_i$ and $c_{ijk}$. For quadratic polynomials, we prove the following main result.

\proclaim{(1.1) Theorem} Let 
$$f(x)=\sum_{1 \leq i < j \leq n} a_{ij}\xi_i \xi_j+ \sum_{i=1}^n b_i \xi_i \quad \text{for} \quad x=\left(\xi_1, \ldots, \xi_n \right).$$
Suppose that for some $0 < \delta < 1$, we have 
$$\sum_{j: \ j \ne i} \left| \Re\thinspace a_{ij} \right| \ \leq \ 1 - \delta, \quad \sum_{j:\ j \ne i} \left| \Im\thinspace a_{ij} \right| \ \leq \ {\delta^2 \over 10} \quad \text{and} \quad \left| \Im\thinspace b_i \right| \ \leq \ 
{\delta^2 \over 10}$$
for $i=1, \ldots, n$.
Then 
$$\sum_{x \in \{-1, 1\}^n} e^{f(x)} \ne 0.$$
\endproclaim

A remark regarding notation: we treat indices $i,j$ in $a_{ij}$ as an unordered pair, so $a_{ij}$ is the coefficient of the monomial $\xi_i \xi_j$ in the quadratic polynomial $f$ and the sum $\sum_{j: \ j \ne i}$ accounts for all 
coefficients $a_{ij}$ of the monomials containing $\xi_i$. For a complex number $z=x+y\sqrt{-1}$, we denote by 
$\Re\thinspace z=x$ and $\Im\thinspace z=y$ the real and imaginary parts of $z$ respectively. 

From Theorem 1.1, the by now standard polynomial interpolation argument, (see \cite{Ba16}, \cite{G+19}, \cite{L+19a}, \cite{L+19b}, \cite{PR17}), produces an algorithm for approximating the partition function $\SS\left(e^f\right)$ when the coefficients $a_{ij}$ and $b_i$ are real and satisfy the condition 
$$\sum_{j: \ j \ne i} |a_{ij}| \ \leq \ 1-\delta \quad \text{for} \quad i=1, \ldots, n, \tag1.1.1$$
where $0 < \delta < 1$ is fixed in advance. As there is no restriction on $b_i$, the sum $\SS\left(e^f\right)$ can be exponentially large in $|b_i|$. To avoid dealing with exponentially large numbers, we assume that we are provided with numbers $e^{b_i}$
for $i=1, \ldots, n$.  Then the  complexity of the algorithm is quasi-polynomial: we approximate $\SS\left(e^f\right)$ within relative error $0 < \epsilon < 1$ in $n^{O(\ln n-\ln \epsilon)}$ time. As Guus Regts explained to the authors, the approach of \cite{PR17} makes our algorithm genuinely polynomial, as opposed to quasi-polynomial, provided the degree of the underlying graph $G$ is fixed in advance.

We describe the algorithm in Section 4 and prove Theorem 1.1 in Section 2. Next, we relate Theorem 1.1 to what is known about the partition function in the Ising model.

\subhead (1.2) The bounds for the zero-free region are asymptotically optimal \endsubhead Let $G$ be a graph with vertices $1, \ldots, n$. For a real number $a$, let us choose 
$$a_{ij}=\cases a &\text{if $\{i,j\}$ is an edge of $G$} \\ 0 & \text{otherwise} \endcases \tag1.2.1$$
and let us choose $b_i=b$ for some $b \in {\Bbb C}$ for all $i=1, \ldots, n$. Let us fix some positive integer $\Delta \geq 3$ and choose either
$$a={1 \over 2} \ln {\Delta \over \Delta -2}$$ 
(all interactions are ferromagnetic)
or 
$$a={1 \over 2} \ln {\Delta-2 \over \Delta}$$
(all interactions are antiferromagnetic).
We consider the partition function $\SS\left(e^f\right)$ as a function of a complex parameter $b$. It is known that for a fixed $\Delta$, as $n$ grows and $G$ ranges over all graphs with the largest degree $\Delta$ of a vertex, the zeros of the univariate function $b \longmapsto \SS\left(e^f\right)$ with either choice of $a$ can get arbitrarily close to $b=0$, see \cite{BM97}, \cite{BG01}, \cite{PR20}, \cite{B+19}. We have
$$\sum_{j: \ j \ne i} \left|\Re\thinspace a_{ij}\right| \ \leq \ {\Delta \over 2} \ln {\Delta \over \Delta-2}.  \tag1.2.2$$
The right hand side approaches $1$ as $\Delta \longrightarrow \infty$, which shows that ``1" in the ``$1-\delta$" bound of Theorem 1.1 cannot be replaced by a larger number.

\subhead (1.3) The bounds for approximation are asymptotically optimal in the antiferromagnetic case \endsubhead As in Section 1.2, let $G$ be a graph of maximum degree $\Delta \geq 3$, let us choose all $b_i=0$ for all $i=1, \ldots, n$ and let us define $a_{ij}$ by (1.2.1).
Suppose that 
$$a \ < \ {1 \over 2} \ln {\Delta -2 \over \Delta}.$$
It is shown in \cite{SS14} and also in \cite{G+16} that the problem of approximating $\SS\left(e^f\right)$ is NP-hard under randomized reduction. Hence unless the computational complexity hierarchy collapses, we cannot approximate 
$\SS\left(e^f\right)$ in quasi-polynomial time in the class of problems where all $b_i=0$ and 
$$\sum_{j: \ j \ne i} |a_{ij}| \ \leq \ 1 + \delta \quad \text{for} \quad i=1, \ldots, n$$
for an arbitrarily small $\delta >0$, fixed in advance.

\subhead (1.4) The ferromagnetic case is special \endsubhead Suppose that $a_{ij} \geq 0$ for all $i, j$ and that $b_i=b$ for $i=1, \ldots, n$ and some complex parameter $b$. Lee and Yang \cite{LY52} showed 
that the zeros of the univariate function $b \longmapsto \SS\left(e^f\right)$ lie on the line $\Re\thinspace b=0$. If the $b_i$ are allowed to vary, then $\SS\left(e^f\right) \ne 0$ as long as $\Re\thinspace b_i > 0$ for $i=1, \ldots, n$. The ferromagnetic case is also special from the complexity point of view: Jerrum and Sinclair \cite{JS93} constructed a randomized polynomial time algorithm approximating $\SS\left(e^f\right)$ when $b_i=b$ for some real $b$ and $i=1, \ldots, n$. Deterministic approximation algorithms of quasi-polynomial (genuinely polynomial, if the maximum degree $\Delta$ of the underlying graph $G$ is fixed in advance) complexity, are constructed in \cite{L+19b}, assuming that $\Re\thinspace b_i \geq \delta > 0$
for some constant $\delta >0$, fixed in advance, and $i=1, \ldots, n$, see also Section 7.4 of \cite{Ba16}. The complexity status of the approximation problem in the ferromagnetic case of $a_{ij} \geq 0$ and $b_i=0$ by a deterministic algorithm appears to be not known.

\subhead (1.5) Correlation decay \endsubhead Theorem 1.1 can be interpreted as saying that there is no phase transition in the Lee - Yang sense \cite{YL52} provided
the parameters $a_{ij}$ and $b_i$ are real and satisfy (1.1.1) for some $0 < \delta < 1$, fixed in advance. There is a related, though not identical, concept of phase transition, based on the disappearance of correlation decay, see \cite{FV18}. If we choose $b_i=0$ for all $i$ and 
define $a_{ij}$ as in (1.2.1), then the correlation decay occurs precisely in the interval 
$$ {1 \over 2} \ln {\Delta-2 \over \Delta} \ < \ a \ < \ {1 \over 2} \ln {\Delta \over \Delta-2},$$
where $\Delta \geq 3$ is the largest degree of a vertex of $G$, see \cite{Z+11}, \cite{S+14}, and \cite{L+19a}. 

We extend the bounds of Theorem 1.1 to cubic polynomials on the Boolean cube $\{-1, 1\}^n$.
\proclaim{(1.6) Theorem}
Let 
$$f(x) =\sum_{1 \leq i < j < k \leq n} c_{ijk} \xi_i \xi_j \xi_k + \sum_{1 \leq i < j \leq n} a_{ij} \xi_i \xi_j + \sum_{i=1}^n b_i \xi_i \quad \text{for} \quad x=\left(\xi_1, \ldots, \xi_n\right).$$
Suppose that for some $0 < \delta < 1/2$, we have 
$$\split &\sum_{j, k:\ j, k \ne i} \left|\Re\thinspace c_{ijk} \right| + \sum_{j: \ j \ne i} \left|\Re\thinspace a_{ij}\right| \ \leq \ 1-\delta, \quad \sum_{j, k:\ j, k \ne i} \left|\Im\thinspace c_{ijk} \right| + \sum_{j: \ j \ne i} \left|\Im \thinspace a_{ij}\right| \ \leq \ {\delta^2 \over 10} \\
& \text{and} \quad \left|\Im\thinspace b_i \right| \ \leq \ {\delta^2 \over 10} \endsplit $$
for $i=1, \ldots, n$.
Then
$$\sum_{x \in \{-1, 1\}^n} e^{f(x)} \ne 0.$$
\endproclaim
Similarly to the quadratic case, we treat indices $i, j, k$ in $c_{ijk}$ as unordered triples of distinct numbers, so that $c_{ijk}$ is the coefficient of the monomial $\xi_i \xi_j \xi_k$ in the cubic polynomial $f$ and the sum 
$\sum_{j, k:\ j, k \ne i}$ accounts for all coefficients $c_{ijk}$ of the monomials containing $\xi_i$. As in the case of quadratic polynomials, we obtain an algorithm of quasi-polynomial $n^{O(\ln n - \ln \epsilon)}$ complexity approximating 
$$\SS\left(e^f\right) = \sum_{x \in \{-1, 1\}^n} e^{f(x)}$$
within relative error $0 < \epsilon < 1$, provided the coefficients $c_{ijk}$, $a_{ij}$ and $b_i$ are real and satisfy 
$$\sum_{j, k:\ j, k \ne i} \left|c_{ijk} \right| + \sum_{j: \ j \ne i} \left|a_{ij}\right| \ \leq \ 1-\delta \quad \text{for} \quad i=1, \ldots, n. \tag1.6.1$$
The bounds here are asymptotically optimal, for the trivial reason that they are asymptotically optimal in the quadratic case when $c_{ijk}=0$ for all $1 \leq i < j < k \leq n$. 

Although Theorem 1.1 is a particular case of Theorem 1.6, we first prove Theorem 1.1 in Section 2 and then describe in Section 3 how to adjust the proof of Theorem 1.1 to obtain Theorem 1.6. We think that this way it is easier to see the idea of the proofs.  In Section 3, we also mention where our method falls short extending the bounds to polynomials of degree 4 and higher.
We don't discuss the corresponding extension of the approximation algorithm of Section 4 from quadratic to cubic polynomials, as it is quite straightforward.

\subhead (1.7) What's new \endsubhead  In the case of a quadratic polynomial, the bound (1.1.1) states that the Lipschitz constant of the non-linear part of $f: \{-1, 1\}^n \longrightarrow {\Bbb R}$ does not exceed $1-\delta$ with respect to the $\ell^1$ metric (which is twice the Hamming metric) on the Boolean cube $\{-1, 1\}^n$. This is where our approach to approximation differs from those of \cite{G+19}, \cite{L+12},  \cite{L+19a}, \cite{L+19b}, \cite{S+14} and \cite{Z+11}, which require a uniform bound on the strength of individual interactions $a_{ij}$, such as
$$\max_{i,j} |a_{ij}| \ \leq \ {1 \over 2} \ln {\Delta \over \Delta-2} \tag1.7.1$$
where $\Delta$ is the maximum degree of the underlying graph.
As we remarked above, our results are asymptotically optimal, when $\Delta \longrightarrow \infty$. Generally, the condition (1.1.1) appears to be more robust than (1.7.1), as (1.1.1) is independent of the degree $\Delta$ and allows individual coefficients 
$a_{ij}$ to be relatively large, as long as the sum for all interactions of any given vertex remains appropriately bounded. Of course, for any particular $\Delta$, the conditions (1.1.1) and (1.7.1) are in general position, as it is easy to construct examples where one holds and the other is violated.

A novelty of our approach with respect to locating zero-free regions of $\SS\left(e^f\right)$, compared to those of \cite{G+19}, \cite{L+19a},  \cite{B+19} and \cite{PR20}, is that we allow all parameters $a_{ij}$ and $b_i$ to vary:
this concerns both Lee - Yang zeros \cite{LY52}, \cite{B+19}, \cite{PR20} of $\SS\left(e^f\right)$ as a function of $b_i$ with fixed $a_{ij}$ and Fisher zeros \cite{L+19a} of $\SS\left(e^f\right)$ as a function of $a_{ij}$ with $b_i$ fixed. It appears that Theorem 1.1 is the first result establishing an asymptotically optimal zero-free region when the interactions $a_{ij}$ are allowed to differ for different pairs $\{i, j\}$ and even to be of different signs, so that we have a mixture of ferromagnetic and antiferromagnetic interactions.

There is much less known in the case of cubic and higher degree polynomials on the Boolean cube, see though \cite{L+16} and \cite{L+19b}. In particular, Theorem 1.6 might be the only result that provides asymptotically optimal bounds for cubic polynomials, though the optimality follows by a trivial reduction to the quadratic case. We note that in contrast to (1.1.1), the left hand side of (1.6.1) is not the Lipschitz constant of the non-linear part of $f$ but an upper bound on that Lipschitz constant.

For polynomials $f$ of degree $d \geq 4$ on the Boolean cube, some, apparently non-optimal, estimates can be found in \cite{Ba17}. Roughly, estimates of \cite{Ba17} require the Lipschitz constant of $f$ to be 
$O(1/\sqrt{d})$. It would be interesting to find out if one can get a bound uniform on $d$.

\head 2. Proof of Theorem 1.1 \endhead

\subhead (2.1) Definitions and notation \endsubhead We consider the Boolean cube $\{-1, 1\}^n$ of vectors $x=\left(\xi_1, \ldots, \xi_n \right)$, where $\xi_i = \pm 1$ for $i=1, \ldots, n$. 
Let us choose a set $I \subset \{1, \ldots, n\}$ and numbers $\sigma_i \in \{-1, 1\}$ for $i \in I$. The set 
$$F=\bigl\{x \in \{-1, 1\}^n: \ \xi_i =\sigma_i \quad \text{for} \quad i \in I \bigr\}$$
is called a {\it face} of $\{-1, 1\}^n$. Any index $i \in \{1, \ldots, n\}\setminus I$ is called a {\it free index} of the face $F$ and the number $n-|I|$ of free indices is called the {\it dimension} of $F$ and denoted $\dim F$.
Indexes $i \in I$ are called {\it fixed}.
For example, if $I =\emptyset$ then the face is the whole cube $\{-1, 1\}^n$ and if $I=\{1, \ldots, n\}$ then the face consists of a single point. Generally, a face of dimension $k$ consists of $2^k$ points.

Let $g: \{-1, 1\}^n \longrightarrow {\Bbb C}$ be a function and let $F \subset \{-1, 1\}^n$ be a face. We define the {\it partial sum of $g$ on $F$} by 
$$\SS (g| F) =  \sum_{x \in F} g(x).$$ We will use the following straightforward identity. Let $F$ be a face of $\{-1, 1\}^n$ of dimension at least $1$. Let $i$ be a free index of $F$ and let 
$F^+ \subset F$ and $F^- \subset F$ be the faces obtained by fixing the $i$-th coordinate $\xi_i$ of vectors $x \in F$ to $1$ and $-1$ respectively. Hence $F^+ \cup F^-=F$, $F^+ \cap F^- =\emptyset$ and
$$\SS(g| F) =  \SS(g | F^+) + \SS(g| F^-). \tag2.1.1$$
Let us fix a real number $0 < \delta < 1$ and a vector $b=\left(b_1, \ldots, b_n\right)$ such that $\left| \Im\thinspace b_i\right| \leq \delta^2/10$ for $i=1, \ldots, n$. We denote by
$\UU(b; \delta)$ the set of all polynomials $f: \{-1, 1\}^n \longrightarrow {\Bbb C}$, 
$$f\left(\xi_1, \ldots, \xi_n \right) = \sum_{1 \leq i < j \leq n} a_{ij} \xi_i \xi_j + \sum_{i=1}^n b_i \xi_i,$$
where 
$$\sum_{j=1}^n \left| \Re\thinspace a_{ij} \right| \ \leq \ 1 - \delta \quad \text{and} \quad \sum_{j=1}^n \left| \Im \thinspace a_{ij} \right| \ \leq \ {\delta^2 \over 10}  \quad \text{for} \quad i=1, \ldots, n.$$
We can view $\UU(b; \delta)$ as a convex subset of ${\Bbb C}^{n \choose 2}$ with a non-empty interior. In particular, for any $f_1, f_2 \in \UU(b; \delta)$ and any $0 \leq \alpha \leq 1$, we have $f \in \UU(b; \delta)$ for
 $f=\alpha f_1 + (1-\alpha) f_2$. Our goal is to show that $\SS\left(e^f \right) \ne 0$ for all $f \in \UU(b; \delta)$.
 
 In what follows, we view non-zero complex numbers $z \in {\Bbb C} \setminus \{0\}$ as vectors in the plane ${\Bbb R}^2={\Bbb C}$ and measure angles between them.
 We prove by induction on the dimension of a face $F \subset \{-1, 1\}^n$ that $\SS\left(e^f | F\right) \ne 0$ and, moreover, if $F^+$ and $F^-$ are faces of $F$ defined above then 
 the angle between complex numbers $\SS\left(e^f | F^+\right) \ne 0$ and $\SS\left(e^f | F^-\right) \ne 0$ is small.

We start with a simple geometric lemma.
\proclaim{(2.2) Lemma} Let $w_+, w_- \in {\Bbb C}\setminus \{0\}$ be non-zero numbers such that the angle between $w_+$ and $w_-$ does not exceed some $0 \leq \theta < \pi$ and let 
$w=w_+ + w_-$. 
\roster 
\item We have $w \ne 0$ and, moreover, 
$${|w_+| + |w_-| \over |w|} \ \leq \ {1 \over \cos(\theta/2)}.$$
\item We have 
$$\left| \Im\thinspace {w_+ - w_- \over w}\right| \ \leq \ \tan {\theta \over 2}.$$
\endroster
\endproclaim
\demo{Proof} Let $v_+$, $v_-$ and $v$ be the orthogonal projections of $w_+$, $w_-$ and $w$ respectively onto the bisector of the angle between $w_+$ and $w_-$. Then 
$$|w| \ \geq \ |v| \ = \left| v_+ + v_-  \right| = |v_+| +  |v_-| \geq \cos(\theta/2) |w_+| + \cos(\theta/2) |w_-|$$
and Part 1 follows.

To prove Part 2, let $u_+ = w_+/w$ and $u_-=w_-/w$. Then 
$$u_+ + u_-=1$$
is real and the angle between $u_+$ and $u_-$ does not exceed $\theta$. Let
$$\alpha = \arg u_+ \quad \text{and} \quad \beta = -\arg u_-.$$
Without loss of generality, we may assume that $0 < \alpha, \beta < \pi/2$ and hence $\alpha + \beta \leq \theta$. Let
$$c= |u_+| + |u_-| \ \leq \ {1 \over \cos (\theta/2)},$$
where the inequality follows from Part 1.
Since $\Im (u_+ + u_-)=0$, we have 
$$(\sin \alpha) |u_+|= (\sin \beta) |u_-| $$
and hence 
$$|u_+|={c  \sin \beta \over \sin \alpha + \sin \beta } \quad \text{and} \quad |u_-| ={c \sin \alpha \over \sin \alpha + \sin \beta}.$$
We have 
$$\left|\Im\thinspace (u_+ - u_-) \right| =\left|\Im \thinspace u_+ \right| + \left| \Im \thinspace u_- \right| = (\sin \alpha) |u_+| + (\sin \beta) |u_-| = {2 c ( \sin \alpha )(\sin \beta) \over \sin \alpha  + \sin \beta}.$$
Since the function $x \longmapsto 1/\sin (x)$ is convex on the interval $0 < x < \pi/2$, the minimum value of 
$${\sin \alpha + \sin \beta \over (\sin \alpha) (\sin \beta)} = {1 \over \sin \alpha} + {1 \over \sin \beta}$$
subject to the constraints $0 < \alpha, \beta < \pi/2$ and $\alpha + \beta \leq \theta$ is attained at $\alpha=\beta = \theta/2$. Consequently, 
$$\left|\Im \thinspace (u_+ -u_-)\right|   \leq  {2 c  \sin^2 (\theta/2) \over 2 \sin (\theta/2)} = c \sin (\theta/2)\ \leq \ \tan {\theta \over 2},$$
which proves Part 2.
{\hfill \hfill \hfill} \qed
\enddemo

Part 1 can be extended to the sum of more than two vectors with pairwise angles not exceeding $\theta$, for which one should require $\theta < 2 \pi/3$, see Lemma  3.6.3 in \cite{Ba16}. 

\proclaim{(2.3) Lemma} Let $F \subset \{-1, 1\}^n$ be a face of $\{-1, 1\}^n$. Suppose that for all $f \in \UU(b; \delta)$ we have 
$\SS \left(e^f| F\right) \ne 0$ and, moreover, the following condition is satisfied: if $i$ is a free index of $F$ and $F^+ \subset F$ and $F^- \subset F$ are the faces obtained by setting 
the $i$-th coordinate $\xi_i=1$ and $\xi_i=-1$ respectively, then the angle between the numbers $\SS\left(e^f| F^+\right) \ne 0$ and $\SS\left(e^f| F^-\right) \ne 0$ does not exceed $\theta$ for
some $0 \leq \theta < \pi$.

Since $\SS \left(e^f | F \right) \ne 0$ for all $f \in \UU(b; \delta)$,
$$f\left(\xi_1, \ldots, \xi_n\right)=\sum_{1 \leq i < j \leq n} a_{ij} \xi_i \xi_j + \sum_{i=1}^n b_i \xi_i,$$
and the set $\UU(b; \delta)$ is simply connected, we can choose a branch of the function $f \longmapsto \ln \SS \left(e^f |F \right)$ for $f \in \UU(b; \delta)$.

Let us fix two indices $1 \leq p < q \leq n$, at most one of which is free for $F$. Then 
$$\left| {\partial \over \partial a_{pq}} \ln \SS \left(e^f | F \right) \right| \ \leq \ {1 \over \cos (\theta/2)} \quad \text{and} \quad 
\left| \Im\thinspace {\partial \over \partial a_{pq}}  \ln \SS \left( e^f | F \right) \right| \ \leq \ \tan {\theta \over 2}$$
for all $f \in \UU(b; \delta)$.
\endproclaim
\demo{Proof}
Differentiating, we get 
$${\partial \over \partial a_{pq}} \ln \SS \left(e^f | F \right) = {\SS \left(\xi_p \xi_q e^f | F \right) \over \SS \left(e^f |F \right)}.$$
Suppose first that neither of the indices $p$ or $q$ is free. Then the value of $\xi_p \xi_q$ is constant for all $x \in F$ and hence 
$${\SS \left(\xi_p \xi_q e^f | F \right) \over \SS \left(e^f |F \right)}= \xi_p \xi_q,$$
from which the conclusions of the lemma follow.

Suppose now that only one of the indices $p$ and $q$, say $q$, is free. Let $F^+ \subset F$ and $F^- \subset F$ be the faces obtained by setting the $q$-th coordinate $\xi_q=1$ and $\xi_q=-1$ respectively.
Applying (2.1.1), we get 
$$ {\SS \left(\xi_p \xi_q e^f | F \right) \over \SS \left(e^f |F \right)}= \xi_p  {\SS \left(\xi_q e^f | F \right) \over \SS \left(e^f |F \right)} = \xi_p  \left( {\SS \left( e^f | F^+ \right) \over \SS \left(e^f |F \right)} - 
{\SS \left( e^f | F^- \right) \over \SS \left(e^f |F \right)} \right).$$
We apply Lemma 2.2 with 
$$w_+=\SS \left( e^f | F^+ \right) \text{and} \quad w_-=\SS \left( e^f | F^- \right).$$
By (2.1.1), we have 
$$w=w_+ + w_- =\SS\left(e^f |F \right)$$
and the angle between $w_+ \ne 0$ and $w_- \ne 0$ does not exceed $\theta$ by the assumption of the lemma.
Applying Part 1 of Lemma 2.2, we get
$$\left| {\partial \over \partial a_{pq}} \ln \SS \left(e^f | F \right) \right| = \left| {w_+ \over w} - {w_- \over w} \right| \ \leq \ {|w_+| + |w_-| \over |w|} \ \leq \
{ 1 \over \cos (\theta/2)}.$$
 Applying Part 2 of Lemma 2.2, we conclude that 
$$\left| \Im\thinspace  {\partial \over \partial a_{pq}}  \ln \SS \left( e^f | F \right) \right| = \left| \Im \thinspace {w_+ -w_- \over w}  \right| \ \leq \ \tan {\theta \over 2}, $$
which completes the proof.
{\hfill \hfill \hfill} \qed
\enddemo

\proclaim{(2.4) Corollary} Let $F \subset \{-1, 1\}^n$ be a face as in Lemma 2.3 and let $f \in \UU(b; \delta)$ be a polynomial. Suppose that $p$ is a fixed index of $F$ and let $g \in \UU(b; \delta)$ be a polynomial obtained by replacing the coefficient $a_{pq}$ in $f$ for some $q$ by $-a_{pq}$. Then the angle between $\SS\left(e^f| F\right) \ne 0$ and $\SS\left(e^g| F \right) \ne 0$ does not exceed
$$\left(2\tan{\theta \over 2} \right) \left| \Re\thinspace a_{pq}\right| + {2\left|\Im\thinspace a_{pq}\right| \over \cos(\theta/2)}.$$
\endproclaim
\demo{Proof} For $-1 \leq t \leq 1$, let $f_t: \{-1, 1\}^n \longrightarrow {\Bbb R}$ be the polynomial obtained by replacing $a_{pq}$ with $t a_{pq}$ in $f$, so $f_1=f$ and $f_{-1}=g$ and 
$f_t \in \UU(b; \delta)$ for $-1 \leq t \leq 1$. Then 
$$\ln \SS\left(e^f | F\right) - \ln \SS \left(e^g | F \right) = \int_{-1}^1 {d \over dt} \ln \SS\left(e^{f_t}|F \right) \ dt =a_{pq} \int_{-1}^1 {\partial \over \partial a_{pq}} \ln \SS\left(e^{f_t}|F \right) \ dt .$$
Consequently,
$$\split &\left| \Im \ln \SS\left(e^f | F\right) -\Im  \ln \SS \left(e^g | F \right) \right| \ \leq \ |\Re\thinspace a_{pq}| \int_{-1}^1 \left| \Im\thinspace  {\partial \over \partial a_{pq}} \ln \SS\left(e^{f_t}|F \right) \right| \ dt\\&\qquad  +
|\Im\thinspace a_{pq}| \int_{-1}^1 \left|\thinspace  {\partial \over \partial a_{pq}} \ln \SS\left(e^{f_t}|F \right) \right| \ dt. \endsplit$$
Applying Lemma 2.3, we complete the proof.
{\hfill \hfill \hfill} \qed
\enddemo

\subhead (2.5) Proof of Theorem 1.1 \endsubhead First, we show that there is $0 < \theta < \pi/2$ such that 
$$ 2  (1-\delta) \tan {\theta \over 2} + {\delta^2 \over 5\cos (\theta/2)}  +{\delta^2 \over 5} \ \leq \ \theta. \tag2.5.1$$
Indeed, we can just choose $\theta=\delta/2$. Using that 
$$\cos x \ \geq \ \cos \left({1 \over 4}\right) \ \geq \ {9 \over 10} \quad \text{for} \quad 0 \leq x \leq {1 \over 4}$$
and that
$$\tan x \ \leq \ x\left(1+{x\over 10} \right) \quad \text{for} \quad 0 \leq x \leq {1 \over 4}$$
(since the function $x^{-2}\left(\tan x -x\right)$ is increasing for $x >0$, it suffices to check the inequality for $x=1/4$),
we obtain
$$\split &2  (1-\delta) \tan {\theta \over 2} + {\delta^2 \over 5\cos (\theta/2)}  +{\delta^2 \over 5} =2  (1-\delta) \tan {\delta \over 4} + {\delta^2 \over 5\cos (\delta/4)} +{\delta^2 \over 5}\\ \leq \ &2(1-\delta) {\delta \over 4} \left(1+{\delta \over 40}\right) + {19 \delta^2 \over 45}={\delta \over 2}
\left(1-{39\delta \over 40}-{\delta^2 \over 40}\right)+ {19 \delta^2 \over 45} \ \leq \ {\delta \over 2}=\theta.  \endsplit$$
We prove by induction for $k=0, \ldots, n$ the following statement.
\bigskip
Let $F \subset \{-1, 1\}^n$ be a face of dimension $k$. Then $\SS\left(e^f | F \right) \ne 0$. Moreover, if $k >0$ and $q$ is a free index of $F$ then the following holds. Let $F^+ \subset F$ and $F^- \subset F$ be the faces of $F$ 
obtained by fixing the $q$-th variable $\xi_q$ to $\xi_q=1$ and $\xi_q=-1$ respectively. Then the angle between $\SS\left(e^f |F^+\right) \ne 0$ and $\SS\left( e^f |F^-\right)\ne 0$ does not exceed $\theta$.
\bigskip
The statement clearly holds for $k=0$. Suppose that $k >1$, let $q$ be a free index of $F$ and let $F^+, F^- \subset F$ be the corresponding faces, so $\dim F^+ = \dim F^-=k-1$. By the induction hypothesis, we have 
$\SS\left(e^f | F^+\right) \ne 0$ and $\SS\left(e^f | F^-\right) \ne 0$. Moreover, $\SS\left(e^f | F^-\right)=\SS\left(e^g| F^+\right)$, where $g$ is obtained from 
$$f\left(\xi_1, \ldots, \xi_n\right)=\sum_{1 \leq i < j \leq n} a_{ij} \xi_i \xi_j + \sum_{i=1}^n b_i \xi_i$$
by replacing the coefficients $a_{pq}$  for $p \ne q$ by $-a_{pq}$ and the coefficient $b_q$ by $-b_q$. Since $q$ is a fixed index for $F^+$, replacing $b_q$ by $-b_q$ leads to multiplying $\SS\left(e^f|F^+\right)$ by $e^{-2b_q \xi_q}$, which results in the rotation of $\SS\left(e^f| F^{+}\right)$ by an angle of $2(\Im\thinspace b_q ) \leq \delta^2/5$ and some scaling. By Corollary 2.4, replacing all $a_{pq}$  with $-a_{pq}$ leads to a rotation of $\SS\left(e^f | F^+\right)$ by at most an angle of 
$$ \left(2 \tan {\theta \over 2} \right) \sum_{p: \ p\ne q} |\Re \thinspace a_{pq} | + {2 \over \cos(\theta/2)} \sum_{p:\ p \ne q} |\Im \thinspace a_{pq}| \ \leq \ 2(1-\delta) \tan {\theta \over 2} +{\delta^2 \over 5 \cos(\theta/2)}$$
and also some scaling.
Hence the angle between $\SS\left(e^f |F^+\right) \ne 0$ and $\SS\left(e^f | F^-\right) =\SS\left(e^g| F^+\right) \ne 0$ does not exceed 
$$2(1-\delta) \tan {\theta \over 2} +{\delta^2 \over 5 \cos(\theta/2)} + {\delta^2 \over 5} \ \leq \ \theta$$
by (2.5.1).
Since
$$\SS\left(e^f | F\right) = \SS\left(e^f| F^+\right) + \SS\left(e^f | F^- \right),$$
applying Part 1 of Lemma 2.2 with 
$$w_+=\SS\left(e^f | F^+ \right) \quad \text{and} \quad w_-=\SS\left(e^f | F^- \right),$$
we conclude that 
$$\SS\left(e^f | F \right) \ne 0.$$
This concludes the proof of the induction step and hence of Theorem 1.1.
{\hfill \hfill \hfill} \qed

\head 3. Proof of Theorem 1.6 \endhead 

First, we obtain a version of Lemma 2.2.

\proclaim{(3.1) Lemma} Let $v_{++}, v_{+-}, v_{-+}, v_{--} \in {\Bbb C} \setminus \{0\}$ be non-zero numbers, such that the angles 
$$\split &\text{between} \quad v_{++} \quad \text{and} \quad v_{+-}, \\
&\text{between} \quad v_{++} \quad \text{and} \quad v_{-+}, \\ 
&\text{between} \quad v_{--} \quad \text{and} \quad v_{+-}, \\
&\text{between} \quad v_{--} \quad \text{and} \quad v_{-+} \endsplit$$
do not exceed some $0 \leq \theta < \pi/2$ and let $v=v_{++}+ v_{+-}+v_{-+}+v_{--}$.
\roster
\item We have $v \ne 0$ and, moreover, 
$${|v_{++}|+|v_{+-}|+|v_{-+}|+|v_{--}| \over |v|} \ \leq \ {1 \over \cos \theta}.$$
\item We have 
$$\left|\Im\thinspace {v_{++}-v_{+-}-v_{-+}+v_{--} \over v} \right| \ \leq \ \tan {\theta \over 2}.$$
\endroster
\endproclaim
\demo{Proof} Clearly, the angle between any two vectors of $v_{++}, v_{+-}, v_{-+}$ and $v_{--}$ does not exceed $2 \theta < \pi$. Let $K \subset {\Bbb R}^2 = {\Bbb C}$ be the angle (cone) spanned by the four vectors. 
Hence the angular measure of $K$ is at most $2\theta$. 

The proof of Part 1 goes exactly as in the proof of Lemma 2.2. \newline
We let $u_{++}, u_{+-}, u_{-+}, u_{--}$ and $u$ be the orthogonal projections of the vectors $v_{++}, v_{+-}, v_{-+}, v_{--}$ and $v$ respectively onto the bisector of $K$. Then 
$$|v| \ \geq \ |u|= |u_{++}| + |u_{+-}|+|u_{-+}| + |u_{--}| \ \geq \ \cos(\theta) \left(|v_{++}| + |v_{+-}| +|v_{-+}| + |v_{--}| \right)$$
and the proof of Part 1 follows.

The proof of Part 2 is obtained by the application of Part 2 of Lemma 2.2. Namely, let 
$$w_+ = v_{++} + v_{--} \quad \text{and} \quad w_-=v_{+-} + v_{-+},$$
so that 
$${v_{++}-v_{+-}-v_{-+}+v_{--} \over v} = {w_+ - w_- \over w_+ + w_-}.$$
We claim that the angle between $w_+$ and $w_-$ does not exceed $\theta$. Indeed, without loss of generality, we assume that one side of the angle $K$ is spanned by $v_{++}$. If the other side of $K$ is spanned by one of the vectors $v_{+-}$ or $v_{-+}$ or the angle is just the ray spanned by $v_{++}$, then the angular measure of $K$ does not exceed $\theta$ and hence the angle between $w_+$ and $w_-$ also does not exceed $\theta$. Suppose now that the other side of $K$ is spanned by $v_{--}$. Without loss of generality, we assume that listed in the clockwise order, the vectors in $K$ are $v_{++}, v_{+-}, v_{-+}$ and $v_{--}$. 
We consider the following two angles: angle $K_1$ spanned by $v_{++}$ and $v_{-+}$ and angle $K_2$ spanned by $v_{+-}$ and $v_{--}$. We note that the angular measure of $K_1$ does not exceed $\theta$ and the angular measure of $K_2$ does not exceed $\theta$. Vector $w_-$ lies in the angle spanned by $v_{+-}$ and $v_{-+}$ and hence in $K_1 \cap K_2$.  On the other hand, since 
$K =K_1 \cup K_2$, vector $w_+$ lies in at least one of the angles $K_1$ or $K_2$. In either case, the angle between $w_+$ and 
$w_-$ does not exceed $\theta$. Applying Part 2 of Lemma 2.2, we conclude that 
$$\left| \Im\thinspace {v_{++}-v_{+-}-v_{-+}+v_{--} \over v}\right| = \left|\Im\thinspace {w_+ - w_- \over w_+ + w_-}\right| \ \leq \ \tan {\theta \over 2},$$
as claimed.
{\hfill \hfill \hfill} \qed
\enddemo

Similar to Section 2, we introduce a set $\UU(b; \delta) \subset {\Bbb C}^{{n \choose 3} +{n \choose 2}}$ of cubic polynomials
$$f(x)=\sum_{1 \leq i < j < k \leq n} c_{ijk} \xi_i \xi_j \xi_k + \sum_{1 \leq i < j \leq n} a_{ij} \xi_i \xi_j + \sum_{i=1}^n b_i \xi_i \quad \text{for} \quad x=\left(\xi_1, \ldots, \xi_n\right)$$
such that 
$$\split &\sum_{j,k:\ j, k \ne i} \left|\Re\thinspace c_{ijk}\right| + \sum_{j:\ j \ne i} \left|\Re\thinspace a_{ij}\right| \ \leq \ 1-\delta, \\
& \sum_{j,k:\ j, k \ne i} \left|\Im\thinspace c_{ijk}\right| + \sum_{j:\ j \ne i} \left|\Im\thinspace a_{ij}\right| \ \leq \ {\delta^2 \over 10} \quad \text{and} \\
&|\Im\thinspace b_i| \leq {\delta^2 \over 10} \quad \text{for} \quad i=1, \ldots, n. \endsplit$$

Next, we need a refinement of Lemma 2.3.
\proclaim{(3.2) Lemma} Let $F \subset \{-1, 1\}^n$ be a face of $\{-1, 1\}^n$. Suppose that for all $f \in \UU(b; \delta)$, we have $\SS\left(e^f| F\right) \ne 0$ and, moreover, the following condition is satisfied: 
if $G \subseteq F$ is a face, if $i$ is a free index of $G$ and $G^+ \subset G$ and $G^- \subset G$ are the faces obtained by setting the $i$-th coordinate $\xi_i=1$ and $\xi_i=-1$ respectively, then the angle between
numbers $\SS\left(e^f |G^+\right) \ne 0$ and $\SS\left(e^f|G^-\right)\ne 0$ does not exceed $\theta$ for some $0 \leq \theta < \pi/2$. We pick a branch of $\ln \SS\left(e^f | F\right)$ for $f \in \UU(b; \delta)$.
\roster
\item
Let us fix three indices $1 \leq p < q < r \leq n$, at most two of which are free for $F$. Then 
$$\left| {\partial \over \partial c_{pqr}} \ln \SS\left(e^f | F \right) \right| \ \leq \ {1 \over \cos \theta} \quad \text{and} \quad \left| \Im\thinspace {\partial \over \partial c_{pqr}} \ln \SS\left(e^f | F \right) \right| \ \leq \ \tan {\theta \over 2}.$$
\item Let us fix two indices $1 \leq p < q \leq n$, at most one of which is free for $F$. Then 
$$\left| {\partial \over \partial a_{pq}} \ln \SS\left(e^f | F \right) \right| \ \leq \ {1 \over \cos (\theta/2)} \ \leq \ {1 \over \cos \theta} \quad \text{and} \quad \left| \Im\thinspace {\partial \over \partial a_{pq}} \ln \SS\left(e^f | F \right) \right| \ \leq \ \tan {\theta \over 2}.$$
\endroster
\endproclaim
\demo{Proof} If there is at most one free index, the proof follows that of Lemma 2.3. Suppose that exactly two indices, say $p$ and $q$ are free, while $r$ is fixed.
Then 
$${\partial \over \partial c_{pqr}} \ln \SS\left(e^f | F \right) ={\SS\left(\xi_p \xi_q \xi_r e^f | F\right) \over \SS\left(e^f | F \right)}= {\xi_r \SS\left(\xi_p \xi_q  e^f | F\right) \over \SS\left(e^f | F \right)}.$$
Let us define four faces $F^{++}, F^{+-}, F^{-+}, F^{--} \subset F$ as follows:
$$\split &F^{++} \subset F\quad  \text{is obtained from $F$ by setting} \quad \xi_p=1, \xi_q=1; \\ &F^{+-} \subset F \quad \text{is obtained from $F$ by setting} \quad  \xi_p=1, \xi_q=-1; \\
&F^{-+} \subset F \quad \text{is obtained from $F$ by setting } \quad  \xi_p=-1, \xi_q=1 \\ &F^{--} \subset  F \quad \text{is obtained from $F$ by setting} \quad  \xi_p=-1, \xi_q=-1.
  \endsplit$$
We let 
$$\split &v_{++} = \SS\left(e^f | F^{++}\right), \\ &v_{+-}=\SS\left(e^f| F^{+-} \right), \\ &v_{-+}=\SS\left(e^f| F^{-+}\right), \\ &v_{--}=\SS\left(e^f| F^{--}\right), \endsplit$$
so that 
$$\SS\left(e^f |F \right)=v_{++} +v_{+-}+v_{-+}+v_{--}$$
and
$${\partial \over \partial c_{pqr}} \ln \SS\left(e^f | F \right)=\xi_r {v_{++} - v_{+-} -v_{-+}+v_{--} \over v_{++} + v_{+-} +v_{-+} + v_{--}}.$$
The proof now follows by Lemma 3.1.
{\hfill \hfill \hfill} \qed
\enddemo

As in Section 2, we obtain a corollary.
\proclaim{(3.3) Corollary} Let $F \subset \{-1, 1\}^n$ be a face as in Lemma 3.2 and let $f \in \UU(b; \delta)$ be a polynomial. 
\roster
\item Suppose that $r$ is a fixed index of $F$ and let $g \in \UU(b; \delta)$ be a polynomial obtained from $f$ by replacing the coefficient $c_{pqr}$ in $f$ for some $p$ and $q$ by $-c_{pqr}$. Then 
the angle between $\SS\left(e^f | F\right) \ne 0$ and $\SS\left(e^g | F \right) \ne 0$ does not exceed 
$$\left(2 \tan {\theta \over 2} \right) \left| \Re\thinspace c_{pqr}\right| + {2 \left| \Im\thinspace c_{pqr}\right| \over \cos \theta}.$$
\item Suppose that $p$ is a fixed index of $F$ and let $g \in \UU(b; \delta)$ be a polynomial obtained from $f$ by replacing the coefficient $a_{pq}$ in $f$ for some $q$ by $-a_{pq}$. 
Then the angle between $\SS\left(e^f | F\right) \ne 0$ and $\SS\left(e^g | F \right) \ne 0$ does not exceed 
$$\left(2 \tan {\theta \over 2} \right) \left| \Re\thinspace a_{pq}\right| + {2 \left| \Im\thinspace a_{pq}\right| \over \cos \theta}.$$
\endroster
\endproclaim
\demo{Proof} The proof is similar to that of Corollary 2.4, only that we use Lemma 3.2 instead of Lemma 2.3.
{\hfill \hfill \hfill} \qed 
\enddemo

\subhead (3.4) Proof of Theorem 1.6 \endsubhead First, we show that there is $0 \leq \theta < \pi/4$ such that 
$$2(1-\delta) \tan{\theta \over 2} + {\delta^2 \over 5\cos \theta} + {\delta^2 \over 5} \ \leq \ \theta.$$
As in Section 2.5, we can choose $\theta=\delta/2$ (note that we assumed that $\delta < 1/2$ and hence we can use the estimates of Section 2.5 verbatim). 

We prove by induction the same statement as in Section 2.5. The only difference is that if $r$ is a free index of $F$ and $F^+, F^- \subset F$ are the faces obtained by setting $\xi_r=1$ and $\xi_r=-1$ respectively,
then $\SS\left(e^f | F^+\right)=\SS\left(e^g |F^-\right)$, where $g$ is obtained from $f$ by changing all coefficients $c_{pqr}$ to $-c_{pqr}$, all coefficients $a_{qr}$ to $-a_{qr}$ and the coefficient $b_r$ to $-b_r$.
We then use Corollary 3.3 instead of Corollary 2.4.
{\hfill \hfill \hfill} \qed

\subhead (3.5) Polynomials of higher degree \endsubhead We don't know if the straightforward extensions of Theorems 1.1 and 1.6 hold for polynomials $f: \{-1, 1\}^n \longrightarrow {\Bbb C}$ of degree 4 and higher
(for some non-optimal estimates in that case, see \cite{Ba17}). The reason our argument does not extend to the case of $\deg f \geq 4$ is that we are lacking the appropriate version of Part 2 of Lemma 3.1. In the case of
$\deg f=4$, that version would have looked as follows: given eight non-zero complex numbers $$v_{+++},\ v_{++-},\ v_{+-+},\ v_{+--},\ v_{-++},\ v_{--+},\ v_{-+-},\ v_{---},$$ such that the angle between any two of them whose indices differ in one position does not exceed some small number $\theta >0$, we have
$$\left|\Im\thinspace {v_{+++}-v_{++-}-v_{+-+}-v_{-++} +v_{+--}+v_{-+-}+v_{--+} -v_{---} \over v_{+++}+v_{++-}+v_{+-+}+v_{-++} +v_{+--}+v_{-+-}+v_{--+} +v_{---}}\right| \ {\buildrel ? \over \leq } \tan {\theta \over 2}.$$
However, the above inequality cannot hold. One can easily construct a counterexample by choosing $|v_{+++}|=|v_{---}|=1$, the angle between $v_{+++}$ and $v_{---}$ equal $3\theta$, while choosing all other vectors to have a very small length. In that case, the left hand side should be about $\tan (3\theta/2)$ and our proof cannot proceed unless we tighten the bound for the sum of the absolute values of the coefficients.This failure of the approach does not, of course, exclude a possibility that the extensions of Theorems 1.1 and 1.6 still hold for higher degree polynomials.

\head 4. Approximation \endhead

Let
$$f(x)=\sum_{1 \leq i < j \leq n} a_{ij} \xi_i \xi_j +\sum_{i=1}^n b_i \xi_i \quad \text{for} \quad x=\left(\xi_1, \ldots, \xi_n\right).$$
We suppose that the coefficients $a_{ij}$ and $b_i$ are real and satisfy 
$$\sum_{j:\ j \ne i} |a_{ij}| \ \leq \ 1-\delta \quad \text{for} \quad i=1, \ldots, n$$
and some $0 < \delta < 1$, fixed in advance (there are no restrictions on $b_i$). Here we sketch an algorithm for approximating 
$$\SS\left(e^f\right)=\sum_{x \in \{-1, 1\}^n} e^{f(x)}.$$
Without loss of generality, we assume that $\delta \leq 0.1$
and that $n$ is large enough,
$$n \ \geq \ {100 \over \delta^2},$$
since for smaller $n$ the sum $\SS\left(e^f\right)$ can be computed by brute force.

We use the method of polynomial interpolation. The gist of the method is summarized in the following lemma, see Section 2.2 of \cite{Ba16} and \cite{Ba19}.

\proclaim{(4.0) Lemma} Let $U \subset {\Bbb C}$ be a connected open set containing $0$ and $1$. Then there is a constant $\gamma(U) > 0$ such that the following holds.
Suppose that
$$p(z)=\sum_{k=0}^n a_k z^k, \quad n \geq 2,$$
is a polynomial such that $p(z) \ne 0$ for all $z \in {\Bbb C}$. Then, for every $0 < \epsilon < 1$, the value of $p(1)$, up to relative error $\epsilon$, is determined by the 
coefficients $a_k$ with 
$$k \ \leq \ \gamma(U) (\ln n - \ln \epsilon)$$
and, moreover, can be computed from those coefficients in $n^{O(1)}$ time.
\endproclaim 

Here we say that two complex numbers $w_1 \ne 0$ and $w_2 \ne 0$ approximate each other within relative error $0 < \epsilon < 1$, if we can write 
$w_1=e^{z_1}$ and $w_2=e^{z_2}$ where $|z_1-z_1| \leq \epsilon$. 

In our case, $U \subset {\Bbb C}$ is a sufficiently small neighborhood of the interval $[0, 1] \subset {\Bbb C}$. Our first goal is to write $\SS\left(e^{f}\right)$ as 
$p(1)$ for a univariate polynomial of $p$ with $\deg p=n^{O(1)}$.

\subhead (4.1) Writing $\SS\left(e^f\right)$ as $p(1)$ for a univariate polynomial $p(z)$ \endsubhead  
First, we rewrite $\SS\left(e^f\right)$ as a polynomial in some new variables. 
We have 
$$\split e^{f(x)} =& \left(\prod_{1 \leq i < j \leq n} e^{a_{ij} \xi_i \xi_j}  \right) \left( \prod_{i=1}^n e^{b_i \xi_i} \right) \\= &\exp\left\{ - \sum_{1 \leq i < j \leq n} a_{ij} \right\} 
 \left(\prod_{1 \leq i < j \leq n} e^{a_{ij}(\xi_i \xi_j +1)}\right) \left(\prod_{i=1}^n e^{b_i \xi_i} \right). \endsplit$$
 Let 
 $$c_{ij} = e^{a_{ij}/n^2}-1 \quad \text{for} \quad 1 \leq i < j \leq n. \tag4.1.1$$
 For given $b_1, \ldots, b_n$, we consider a polynomial
 $$P_b\left(c\right) = \sum_{\xi_1, \ldots, \xi_n =\pm 1}  \left( \prod_{1 \leq i < j \leq n} \left(c_{ij}+1\right)^{n^2(\xi_i \xi_j +1)} \right)  \left(\prod_{i=1}^n e^{b_i \xi_i} \right)$$
 of degree $n^3(n-1)$ in the vector $c$ of ${n \choose 2}$ variables $c_{ij}$. 
 Thus we have 
 $$\SS\left(e^f\right) =\exp\left\{ - \sum_{1 \leq i < j \leq n} a_{ij} \right\}  P_b\left(c\right)$$ and 
we want to approximate $P_b\left(c\right)$. 

Next, for given $c$ and $b_i$, we consider a univariate polynomial 
$$p(z)=P_b\left(z c\right) =  \sum_{\xi_1, \ldots, \xi_n =\pm 1}  \left( \prod_{1 \leq i < j \leq n} \left(z c_{ij}+1\right)^{n^2(\xi_i \xi_j +1)} \right)  \left(\prod_{i=1}^n e^{b_i \xi_i} \right) \tag4.1.2$$
of a complex variable $z$. Our goal is to approximate $p(1)$.

\subhead (4.2) Showing that $p(z) \ne 0$ in a neighborhood of $[0, 1] \subset {\Bbb C}$ \endsubhead
Our next goal is to show that 
$$p(z) \ne 0 \quad \text{provided} \quad |\Re\thinspace z| \leq 1+\delta^2 \quad \text{and} \quad |\Im\thinspace z| \ \leq \ {\delta^2 \over 80}. \tag4.2.1$$
Since $|a_{ij}| \leq 1$ and $n \geq 10$, from (4.1.1) we deduce that 
$$\left|c_{ij}- {a_{ij} \over n^2} \right| \ \leq \ {1 \over n^4} \quad \text{for all} \quad i, j\tag4.2.2$$
(follows from the Taylor series expansion of $e^{a_{ij}/n^2}$).
In particular, 
$$|c_{ij}| \ \leq \ {2 \over n^2} \quad \text{for all} \quad i, j. \tag4.2.3$$
From (4.1.1), we obtain
$$a_{ij} = n^2 \ln \left(1+ c_{ij}\right) \quad \text{for all} \quad i, j.$$
Let us choose an arbitrary $z \in {\Bbb C}$ such that $|z| \leq 1+\delta^2$ and let us define 
$$\widehat{a}_{ij}=\widehat{a}_{ij}(z) = n^2 \ln \left(1 + z c_{ij}\right) \quad \text{for all} \quad i,j, \tag4.2.4$$
where we choose the principal branch of the logarithm, so that $\widehat{a}_{ij}(0)=0$. Let 
$$\widehat{f}(x)=\sum_{1 \leq i < j \leq n} \widehat{a}_{ij} \xi_i \xi_j + \sum_{i=1}^n b_i \xi_i \quad \text{for} \quad x=\left(\xi_1, \ldots, \xi_n\right).$$
From Section 4.1, we have 
$$p(z) = \exp\left\{ \sum_{1 \leq i < j \leq n} \widehat{a}_{ij}\right\} \SS\left(e^{\widehat{f}}\right). \tag4.2.5$$
Combining (4.2.2)--(4.2.4), we obtain
$$\widehat{a}_{ij} =z a_{ij}  +\eta_{ij} \quad \text{where} \quad |\eta_{ij}| \leq {3 \over n^2} \quad \text{for all} \quad i,j$$
(follows from the Taylor series expansion of $\ln (1+ z c_{ij})$).
Consequently,
$$\sum_{j: \ j \ne i} \left|\Re\thinspace \widehat{a}_{ij} \right| \ \leq \ (1- \delta)(1+\delta^2)  + {3 \over n} \ \leq \ 1-{\delta \over 2}.$$
In addition, if $|\Im\thinspace z| \leq \delta^2/80$, then 
$$\sum_{j: \ j \ne i} \left|\Im\thinspace \widehat{a}_{ij}\right| \ \leq \ {\delta^2 \over 80} + {3 \over n} \ \leq \ {\delta^2 \over 40}.$$
Hence by Theorem 1.1 we have that $\SS\left(e^{\widehat{f}}\right) \ne 0$ and therefore by (4.2.5) we conclude that (4.2.1) holds.

\subhead (4.3) Computing $p^{(k)}(0)$ \endsubhead As discussed in Section 2.2 of \cite{Ba16}, see also \cite{PR17} and \cite{Ba19} for some enhancements, as soon as $p(z) \ne 0$ in some 
neighborhood $U$ of the interval $[0, 1] \subset {\Bbb C}$, to approximate $p(1)$ within relative error $0< \epsilon < 1$, it suffices to compute the derivatives
$p^{(k)}(0)$ for $k=\gamma(\ln \deg p - \ln \epsilon)$ where $\gamma=\gamma(U) >0$. Since in our case $\deg p \leq n^4$, to approximate 
$\SS\left(e^f\right)=p(1)$ within relative error $0 < \epsilon < 1$, it suffices to compute $p^{k}(0)$ for $k =O\bigl(\ln n - \ln \epsilon\bigr)$, where the implied constant in the ``$O$" notation depends on $\delta$ alone.

From (4.1.2), we obtain
$$\split p^{(k)}(0)=&\sum\Sb \left(k_1, \ldots, k_s\right): \\ k_1 + \ldots + k_s =k\endSb
{k \choose k_1, \ldots, k_s} \left(2n^2\right)_{k_1} \cdots \left(2n^2 \right)_{k_s} \sum_{\xi_1, \ldots, \xi_n =\pm 1} \left(\prod_{i=1}^n e^{b_i \xi_i} \right) 
\\&\qquad \qquad \times \sum\Sb \left\{ i_1 < j_1, i_2 < j_2, \ldots, i_s < j_s \right\} \\ \xi_{i_1} \xi_{j_1} = \ldots =\xi_{i_s}\xi_{j_s} =1 \endSb  
c_{i_1 j_1}^{k_1} \cdots c_{i_s j_s}^{k_s}, \endsplit$$
where the outermost sum is taken over all positive integer vectors $\left(k_1, \ldots, k_s\right)$ such that $k_1+ \ldots + k_s=k$, 
 the innermost sum is taken over all $s$-sets of pairs $i_1 < j_1$, $i_2 < j_2$, $\ldots$, $i_s < j_s$, such that the products $\xi_{i_1} \xi_{j_1}, \ldots, \xi_{i_s} \xi_{j_s}$ are all equal 1's,
 and where
 $$\left(2n^2\right)_{k_i} = 2n^2 \left(2n^2 -1 \right) \cdots \left(2n^2 -k_i+1\right).$$
  It is convenient to rewrite the 
 above sum in the graph notation. Let $K_n$ be the complete undirected graph with set $V =\{1, \ldots, n\}$ of vertices and set $E$ of edges. We assign weight $c_{ij}$ to edge $\{i, j\}$. We call a map 
 $\sigma:\ V \longrightarrow \{-1, 1\}$ {\it consistent} on an edge $u=\{i, j\}$ if $\sigma(i)=\sigma(j)$. Then 
 $$\split p^{(k)}(0)=&\sum\Sb \left(k_1, \ldots, k_s\right): \\ k_1 + \ldots + k_s =k\endSb{k \choose k_1, \ldots, k_s} \left(2n^2\right)_{k_1} \cdots \left(2n^2 \right)_{k_s} 
 \sum_{\sigma:\ V \longrightarrow \{-1, 1\}} \left(\prod_{i=1}^n e^{b_i \sigma(i)} \right) \\ &\qquad \times \sum\Sb \{u_1, \ldots, u_s\} \subset E: \\ \sigma \text{\ is consistent on each\ } u_1, \ldots, u_s \endSb 
 c_{u_1}^{k_1} \cdots c_{u_s}^{k_s}\\ 
 =&\sum\Sb \left(k_1, \ldots, k_s\right): \\ k_1 + \ldots + k_s =k\endSb {k \choose k_1, \ldots, k_s} \left(2n^2\right)_{k_1} \cdots \left(2n^2 \right)_{k_s}  \sum_{\{u_1, \ldots, u_s\} \subset E}  c_{u_1}^{k_1} \cdots c_{u_s}^{k_s} \\
 &\qquad \qquad \sum\Sb \sigma:\ V \longrightarrow \{-1, 1\}: \\  \sigma \text{\ is consistent on each\ } u_1, \ldots, u_s \endSb \left(\prod_{i=1}^n e^{b_i \sigma(i)} \right).\endsplit $$
 For a given $s \leq k$, there are not more than $n^{2s} \leq n^{2k}$ sets of $s$ distinct edges $u_1, \ldots, u_s$. Given such a set $u_1, \ldots, u_s$, let $W \subset V$ be the set of vertices of $u_1, \ldots, u_s$. Then
 $|W|\leq 2s$ and there are at most $2^s$ maps $\sigma_0: W \longrightarrow \{-1, 1\}$ that are consistent on each edge $u_1, \ldots, u_s$. Finally, given such a map $\sigma_0$, we have 
 $$\sum\Sb \sigma:\ V \longrightarrow \{-1, 1\}: \\ \sigma|W =\sigma_0 \endSb \prod_{i=1}^n e^{b_i \sigma(i)} = \prod_{i \in W} e^{b_i \sigma_0(i)} \prod_{i \in V \setminus W} \left( e^{b_i} + e^{-b_i}\right).$$
Summarizing,
$$\split p^{(k)}(0)=&\sum\Sb \left(k_1, \ldots, k_s\right): \\ k_1 + \ldots + k_s =k\endSb {k \choose k_1, \ldots, k_s} \left(2n^2\right)_{k_1} \cdots \left(2n^2 \right)_{k_s} \sum_{\{u_1, \ldots, u_s\} \subset E} c_{u_1}^{k_1} \ldots 
c_{u_s}^{k_s} \\ &\qquad \sum_{\sigma_0:\ W \longrightarrow \{-1, 1\}}  \prod_{i \in W} e^{b_i \sigma_0(i)} \prod_{i \in V \setminus W} \left( e^{b_i} + e^{-b_i}\right),
\endsplit$$
where the outermost sum is taken over at most ${2k-1 \choose k-1}$ positive integer vectors $\left(k_1, \ldots, k_s\right)$ satisfying the equation $k_1 + \ldots + k_s=k$, the next sum is taken over at most 
$n^{2k}$ subsets of $s \leq k$ edges and the innermost sum is taken over at most $2^s \leq 2^k$ consistent maps $\sigma_0$ on the set $W$ of vertices of $\{u_1, \ldots, u_s\}$. Hence the complexity of
computing $p^{(k)}(0)$ is $n^{O(k)}$ and since $k=O(\ln n - \ln \epsilon)$, we obtain an algorithm of quasi-polynomial $n^{O(\ln n - \ln \epsilon)}$ complexity to approximate $\SS\left(e^f\right)$. 

Guus Regts explained to the authors that the argument of \cite{PR17} allows one to obtain a genuinely polynomial, as opposed to quasi-polynomial, algorithm if the degree $\Delta$ of the underlying graph $G$ with edges $\{i, j\}$ where $a_{ij} \ne 0$ is bounded from above in advance. This is based on the observation that computing $p^{(k)}(0)$ in that case reduces to weighted enumeration of {\it connected} subgraphs of $G$ with $O(k)$ vertices, which can be accomplished in polynomial time as long as $\Delta$ is fixed and $k=O(\ln n - \ln \epsilon)$.

The extension of the algorithm to cubic polynomials $f$ is straightforward. We use then Theorem 1.6 instead of Theorem 1.1 to ascertain that the auxiliary univariate polynomial $p(z)$ does not have roots in a neighborhood of $[0, 1] \subset {\Bbb C}$.

\head Acknowledgment \endhead

The authors are grateful to Alistair Sinclair for answering questions regarding  correlation decay in the Ising model, to Guus Regts for explaining how the approach of \cite{PR17} makes the algorithm polynomial, 
as opposed to quasi-polynomial, if the degree of the underlying graph is bounded from above in advance, and to the anonymous referees for their careful reading of the paper and catching 
inaccuracies.


\Refs
\widestnumber\key{AAAAA}

\ref\key{BG01}
\by J.C.A. Barata and P.S. Goldbaum
\paper On the distribution and gap structure of Lee-Yang zeros for the Ising model: periodic and aperiodic couplings
\jour Journal of Statistical Physics 
\vol 103 
\yr 2001
\pages  no. 5-6, 857--891
\endref 

\ref\key{BM97}
\by J.C.A. Barata and D.H.U. Marchetti
\paper Griffiths' singularities in diluted Ising models on the Cayley tree
\jour Journal of Statistical Physics 
\vol 88 
\yr 1997
\pages no. 1-2, 231--268
\endref

\ref\key{Ba16}
\by A. Barvinok
\book Combinatorics and Complexity of Partition Functions
\bookinfo  Algorithms and Combinatorics, 30
\publ Springer
\publaddr Cham
\yr 2016
\endref

\ref\key{Ba17}
\by A. Barvinok
\paper Computing the partition function of a polynomial on the Boolean cube
\inbook A Journey through Discrete Mathematics
\pages 135--164
\publ Springer
\publaddr Cham
\yr 2017
\endref

\ref\key{Ba19}
\by A. Barvinok
\paper Approximating real-rooted and stable polynomials, with combinatorial applications
\jour Online Journal of Analytic Combinatorics
\paperinfo Issue 14 $\#$ 8
\year 2019
\pages 13 pp
\endref

\ref\key{B+19}
\by F. Bencs, P. Buys, L. Guerini and H. Peters
\paper Lee - Yang zeros of the antiferromagnetic Ising Model
\paperinfo preprint {\tt arXiv:1907.07479}
\yr 2019
\endref

\ref\key{FV18}
\by S. Friedli and Y. Velenik
\book Statistical Mechanics of Lattice Systems. A concrete mathematical introduction
\publ Cambridge University Press
\publaddr Cambridge
\yr 2018
\endref

\ref\key{G+16}
\by A. Galanis,  D. \v{S}tefankovi\v{c} and E. Vigoda
\paper Inapproximability of the partition function for the antiferromagnetic Ising and hard-core models
\jour Combinatorics, Probability and Computing
\vol 25 
\yr 2016
\pages no. 4, 500--559
\endref


\ref\key{G+19}
\by H. Guo, J. Liu and P. Lu
\paper Zeros of ferromagnetic 2-spin systems
\paperinfo preprint \newline {\tt arXiv:1907.06156}
\yr 2019
\endref

\ref\key{J+19}
\by V. Jain, A. Risteski and F. Koehler
\paper Mean-field approximation, convex hierarchies, and the optimality of correlation rounding: a unified perspective
\inbook  STOC'19--Proceedings of the 51st Annual ACM SIGACT Symposium on Theory of Computing
\pages 1226--1236
\publ ACM
\publaddr New York
\yr  2019
\endref

\ref\key{JS93}
\by M. Jerrum and A. Sinclair
\paper Polynomial-time approximation algorithms for the Ising model
\jour  SIAM Journal on Computing
\vol 22 
\yr 1993
\pages no. 5, 1087--1116
\endref

\ref\key{LY52}
\by T.D. Lee and C.N. Yang
\paper Statistical theory of equations of state and phase transitions. II. Lattice gas and Ising model
\jour Physical Review (2) 
\vol 87
\yr 1952
\pages  410--419
\endref

\ref\key{L+12}
\by L. Li, P. Lu and Y. Yin
\paper Correlation decay up to uniqueness in spin systems
\inbook Proceedings of the Twenty-Fourth Annual ACM-SIAM Symposium on Discrete Algorithms
\pages 67--84
\publ SIAM
\publaddr Philadelphia, PA
\yr 2012
\endref

\ref\key{L+19a}
\by J. Liu, A. Sinclair and P. Srivastava
\paper Fisher zeros and correlation decay in the Ising model
\jour Journal of Mathematical Physics
\vol 60 
\yr 2019
\pages no. 10, 103304, 12 pp
\endref

\ref\key{L+19b}
\by J. Liu, A. Sinclair and P. Srivastava
\paper The Ising partition function: zeros and deterministic approximation
\jour  Journal of  Statistical Physics 
\vol 174 
\yr 2019
\pages no. 2, 287--315
\endref

\ref\key{L+16}
\by P. Lu, K. Yang and C. Zhang
\paper FPTAS for hardcore and Ising models on hypergraphs
\inbook 33rd Symposium on Theoretical Aspects of Computer Science, 
\bookinfo LIPIcs. Leibniz International Proceedings in Informatics
\vol 47 
\publ Schloss Dagstuhl. Leibniz-Zent. Inform.
\publaddr Wadern
\yr  2016
\pages Art. No. 51, 14 pp
\endref

\ref\key{PR17}
\by V. Patel and G. Regts
\paper Deterministic polynomial-time approximation algorithms for partition functions and graph polynomials
\jour SIAM Journal on Computing 
\vol 46 
\yr 2017
\pages no. 6, 1893--1919
\endref

\ref\key{PR20}
\by H. Peters and G. Regts
\paper Location of zeros for the partition function of the Ising model on bounded degree graphs
\jour Journal of the London Mathematical Society
\vol 101
\pages Issue 2, 765--785
\yr 2020
\endref

\ref\key{S+14}
\by A. Sinclair, P. Srivastava and M. Thurley
\paper Approximation algorithms for two-state anti-ferromagnetic spin systems on bounded degree graphs
\jour Journal of Statistical Physics
\vol 155 
\yr 2014
\pages  no. 4, 666--686
\endref

\ref\key{SS14}
\by A. Sly and N. Sun
\paper Counting in two-spin models on $d$-regular graphs
\jour The Annals of Probability
\vol 42 
\yr 2014
\pages no. 6, 2383--2416
\endref

\ref\key{YL52}
\by C.N. Yang and T.D. Lee
\paper Statistical theory of equations of state and phase transitions. I. Theory of condensation
\jour Physical Review (2) 
\vol 87
\yr 1952
\pages 404--409
\endref

\ref\key{Z+11}
\by J. Zhang, H. Liang and F. Bai
\paper Approximating partition functions of the two-state spin system
\jour Information Processing Letters	
\vol 111
\yr 2011
\pages no. 14, 702--710
\endref


\endRefs
\enddocument
\end